\documentclass[11pt]{article}
\usepackage[T1]{fontenc}
\usepackage{lmodern,amsmath,amssymb,amsthm,mathtools,stmaryrd,hyperref}
\usepackage[margin=1in]{geometry}
\newtheorem{theorem}{Theorem}[section]
\newtheorem{lemma}[theorem]{Lemma}
\newtheorem{remark}[theorem]{Remark}
\newcommand{\HS}{\mathrm{HS}}
\newcommand{\TC}{\mathfrak{L}_1}

\newcommand{\dd}{\,\mathrm d}
\newcommand{\norm}[1]{\left\|#1\right\|}

\newcommand{\ip}[2]{\langle #1,#2\rangle}
\title{Tensor Covariation of Hilbert-Space Stochastic Integrals: A Covariance-Space Formulation}
\author{Sascha Holl\\\texttt{sascha.holl@gmail.com}}
\date{Technical report --- November 2026}
\begin{document}
\maketitle
\begin{abstract}
We formulate the tensor-valued covariation identity for continuous Hilbert-space local martingales in terms of weighted integrands and a contraction obtained from the positive joint covariance operator. This avoids applying an unbounded inverse covariance operator to arbitrary vectors and permits integrands defined on the covariance range with its intrinsic Hilbert norm. We give a trace-norm estimate and a pathwise Kunita--Watanabe inequality. The formulation uses classical Hilbert-space stochastic integration and trace-ideal estimates.
\end{abstract}

\section{Notation and joint covariance}
Let $(\Omega,\mathcal A,\mathrm P)$ be a complete probability space, $T>0$, $I=(0,T]$, and $(\mathcal F_t)_{t\in\overline I}$ a complete right-continuous filtration. All spaces are real separable Hilbert spaces. For continuous local martingales $M$ in $U$ and $N$ in $V$ starting at zero, denote their scalar quadratic variations by $[M]$ and $[N]$. The operator covariation $\llbracket  M,N\rrbracket_t\in\TC(U,V)$ is oriented by
\[
 \ip{\llbracket  M,N\rrbracket_t u}{v}
   =[\ip{M}{u},\ip{N}{v}]_t.
\]
Here $\TC$ is the trace class, and $\HS$ denotes Hilbert--Schmidt operators. For continuous local martingales the optional and predictable brackets agree. We use the standard existence of the predictable trace-class covariance density of a continuous Hilbert-valued local martingale (see \cite{metivier,dpz,cherny}).

Let $L=(M,N)$ in $U\oplus V$ and $\lambda=[M]+[N]$. Relative to the measure $\dd\lambda$ the joint covariance density has block form
\begin{equation}\label{eq:block}
 Q=\begin{pmatrix} A & C^*\\ C & B \end{pmatrix}\geq0,
 \qquad \operatorname{tr}A+\operatorname{tr}B=1
\end{equation}
$\dd\lambda$-almost everywhere almost surely. The blocks $A,B,C$ are predictable trace-class-valued processes. In particular $\llbracket  M,N\rrbracket_t=\int_0^t C_s\dd\lambda_s$. The analogous identities hold for the two diagonal blocks. In the covariance notation, $C=\operatorname{Cor}(M,N)$ with respect to the common clock $\lambda$, and $A=(\operatorname{Cov}M)\,\mathrm d[M]/\mathrm d\lambda$, $B=(\operatorname{Cov}N)\,\mathrm d[N]/\mathrm d\lambda$, where the individual covariance densities are taken with respect to the respective scalar brackets. All assertions are understood up to the appropriate bracket-null sets, with localization when necessary.

\begin{lemma}[Normalized cross covariance]\label{lem:block}
For any positive block operator as in \eqref{eq:block}, there is a unique contraction $R:U\to V$ with
\begin{equation}\label{eq:R}
 C=B^{1/2}RA^{1/2},\qquad R=P_BRP_A,
\end{equation}
where $P_A,P_B$ are the orthogonal projections onto the closures of the ranges of $A,B$. It is predictable in the strong operator sense when the blocks are predictable.
\end{lemma}
\begin{proof}
Positivity applied to $(tu,v)$ for real $t$ gives $|\ip{Cu}{v}|^2\leq\ip{Au}{u}\ip{Bv}{v}$. Thus the form on $\operatorname{ran}A^{1/2}\times\operatorname{ran}B^{1/2}$ defined by $(A^{1/2}u,B^{1/2}v)\mapsto\ip{Cu}{v}$ is well-defined and bounded by the product of the two norms. Extend by continuity to their closures, apply the Riesz representation theorem, and extend by zero on the orthogonal complement of the first range. This proves existence, uniqueness and the contraction bound. For predictability, the operators
\[
 R_\varepsilon=(B+\varepsilon I)^{-1/2}C(A+\varepsilon I)^{-1/2}
\]
are predictable by continuous functional calculus for fixed $\varepsilon>0$. By \eqref{eq:R} they equal $f_\varepsilon(B)Rf_\varepsilon(A)$, where $f_\varepsilon(x)=\sqrt{x/(x+\varepsilon)}$. They converge strongly to $R$ as $\varepsilon\downarrow0$. This provides a strongly operator-predictable version and does not require the generally unbounded inverses of $A$ or $B$.
\end{proof}

\section{Covariance-space integrands}
For a positive trace-class operator $A$, write $U_A=\operatorname{ran} A^{1/2}$ with the norm $\norm{A^{1/2}u}_{U_A}=\norm{P_Au}_U$. It is a Hilbert space with a continuous embedding into $U$. The map $A^{1/2}:U\to U_A$ is a coisometry. An integrand $\Phi:U_A\to H$ in the Hilbert--Schmidt class is therefore encoded by $S=\Phi A^{1/2}:U\to H$, which is Hilbert--Schmidt and satisfies $S=SP_A$. We similarly write $T=\Psi B^{1/2}:V\to J$. In the stochastic setting these effective operators are assumed predictable, and local integrability means
\begin{equation}\label{eq:int}
 \int_0^T\big(\norm{S_s}_2^2+\norm{T_s}_2^2\big)\dd\lambda_s<\infty\quad\text{a.s.}
\end{equation}
We use the usual stochastic integral constructed by energy-norm completion of ambient bounded operator integrands. To see that the covariance-space notation causes no additional restriction, put $\Phi^{\varepsilon}=S(A+\varepsilon I)^{-1/2}$. These are bounded operators on $U$ (pointwise), predictable, and $\Phi^{\varepsilon}A^{1/2}=Sf_\varepsilon(A)\to S$ in Hilbert--Schmidt norm, with the norms bounded by $\norm{S}_2$. Consequently they converge in the stochastic energy norm after localization, by dominated convergence. The same holds for $\Psi$. The integrals $I=\Phi\cdot M$ and $J'=\Psi\cdot N$ are the corresponding limits, uniformly in probability on compact intervals. Using the common clock instead of each individual scalar bracket does not change the integral: the corresponding densities differ by the scalar Radon--Nikodym factors, which cancel in the energy norm. This is also immediate for bounded elementary integrands, followed by completion.

\begin{theorem}[Tensor covariation identity]\label{thm:main}
Under \eqref{eq:int}, almost surely for every $t\in[0,T]$,
\begin{equation}\label{eq:main}
 \llbracket  I,J'\rrbracket_t=\int_0^t T_sR_sS_s^*\dd\lambda_s.
\end{equation}
The integral is a pathwise Bochner integral in $\TC(H,J)$. Its density satisfies
\begin{equation}\label{eq:bound}
 \norm{TRS^*}_1\leq\norm{T}_2\norm{S}_2,
\end{equation}
and $\dd[I]=\norm{S}_2^2\dd\lambda$, $\dd[J']=\norm{T}_2^2\dd\lambda$. If the integrands are bounded on the ambient spaces $U,V$, the density reduces to $\Psi C\Phi^*$. The latter expression is not in general a composition of operators on the ambient spaces for intrinsically defined integrands.
\end{theorem}
\begin{proof}
The trace-ideal product bound and Lemma~\ref{lem:block} prove \eqref{eq:bound}. Predictability of the trace-class density follows by approximating $S,T$ in Hilbert--Schmidt norm by finite-dimensional coordinate truncations and using strong operator predictability of $R$. This also proves pathwise Bochner integrability by Cauchy--Schwarz and \eqref{eq:int}.

For elementary bounded integrands on a common deterministic partition, scalarizing against $h\in H$ and $j\in J$ gives the usual scalar bracket rule, hence
\[
 [\ip{\Phi\cdot M}{h},\ip{\Psi\cdot N}{j}]_t
 =\int_0^t\ip{C_s\Phi_s^*h}{\Psi_s^*j}\dd\lambda_s.
\]
Indeed, on each interval the integrands are measurable at the left endpoint; one can first approximate their adjoint evaluations by finite-valued random vectors and then use the scalar bracket identity. Since $C=B^{1/2}RA^{1/2}$, the right side is the scalar pairing of \eqref{eq:main}.

The identity extends by stochastic energy-norm approximation. Explicitly, if $S_n,T_n$ are effective elementary integrands converging in the energy norms after localization, then
\begin{align*}
 &\int_0^T\norm{T_nRS_n^*-TRS^*}_1\dd\lambda\\
 &\quad\leq \norm{T_n-T}_{L^2(\lambda;\HS)}\norm{S_n}_{L^2(\lambda;\HS)}
 +\norm{T}_{L^2(\lambda;\HS)}\norm{S_n-S}_{L^2(\lambda;\HS)}.
\end{align*}
This converges to zero in probability under energy convergence. On the other hand, the stochastic integrals converge uniformly in probability and their differences have scalar brackets converging to zero in probability. The scalar Kunita--Watanabe inequality applied to these differences then gives convergence of every paired covariation uniformly in probability. The energy completion and the regularized approximation above include all integrands in \eqref{eq:int}; equivalently one can first stop when the integrated energies and the clock exceed fixed bounds and use the square-integrable theory there. Hence for each $h,j$ the scalar pairings of \eqref{eq:main} agree indistinguishably. Countable dense subsets and continuity identify the trace-class operators for all $t$ simultaneously. Finally the same argument with identical martingales and integrands, or the usual diagonal covariance formula of the stochastic integral, gives the scalar brackets by taking the trace of $SS^*$ and $TT^*$.
\end{proof}
\begin{remark}
The formula includes singular covariance, nonclosed covariance ranges and unequal noise spaces. It avoids expressions such as $B C A^{-1}$, whose domain may not contain the adjoints of covariance-space integrands. The contraction $R$ is dimensionless and is defined using the joint covariance, rather than the total variation of the cross covariation alone.
\end{remark}

\section{Pathwise Kunita--Watanabe inequality}
\begin{theorem}\label{thm:kw}
Let $Z=\llbracket  M,N\rrbracket$ and $|\dd Z|$ be the total variation of its trace-class-valued Stieltjes measure. Outside a single exceptional set, for every pair of Borel functions $f,g$ on $[0,T]$ with finite right-hand side,
\[
 \int |fg|\,|\dd Z|
 \leq\left(\int |f|^2\dd[M]\right)^{1/2}
      \left(\int |g|^2\dd[N]\right)^{1/2}.
\]
\end{theorem}
\begin{proof}
Apply \eqref{eq:R} and the trace-ideal product inequality to obtain $\norm{C}_1\leq\sqrt{\operatorname{tr}A}\sqrt{\operatorname{tr}B}$. Since $\dd Z=C\dd\lambda$, its variation has density $\norm{C}_1$ with respect to $\dd\lambda$. Cauchy--Schwarz proves the result, using $\dd[M]=\operatorname{tr}A\dd\lambda$ and the analogous identity for $N$. The argument is pathwise and simultaneously valid for all such test functions.
\end{proof}
\begin{remark}[Deterministic version without the Radon--Nikod\'ym property]
Let $\mu,\nu$ be finite positive measures, and $m$ a Banach-valued countably additive measure of bounded variation with $\norm{m(E)}\leq\sqrt{\mu(E)\nu(E)}$ for every measurable $E$. For a finite measurable partition $(E_i)$ of $E$ and any $r>0$,
\[
 \sum_i\norm{m(E_i)}\leq\tfrac12\big(r\mu(E)+r^{-1}\nu(E)\big).
\]
Taking the supremum over partitions gives $|m|\leq(r\mu+r^{-1}\nu)/2$ as measures. With $\eta=\mu+\nu$, scalar Radon--Nikod\'ym derivatives and the countability of positive rational $r$ imply $\dd|m|/\dd\eta\leq\sqrt{(\dd\mu/\dd\eta)(\dd\nu/\dd\eta)}$ almost everywhere, by minimizing over $r$. The same weighted inequality follows. No Banach-valued Radon--Nikod\'ym derivative, and thus no Radon--Nikod\'ym property of the range space, is needed. The hypothesis for all measurable sets follows from the interval hypothesis in the original notes by the monotone class theorem: first extend to disjoint finite unions of half-open intervals by Cauchy--Schwarz, then to the generated sigma algebra using continuity of finite measures and vector measures of bounded variation.
\end{remark}

\section{Scope}
The common-clock representation provides an operator-domain-safe formulation of the tensor covariation identity for covariance-space integrands. The underlying Hilbert-space martingale integration and Kunita--Watanabe theory are classical; see, in particular, Asperti \cite{asperti} and Rehmeier \cite{cherny}. No assertion of priority is made.

\end{document}